\documentclass[12pt]{amsart}

\usepackage{tikz}
\usetikzlibrary{calc}
\usepackage{ae} 
\usepackage[T1]{fontenc}
\usepackage[cp1250]{inputenc}
\usepackage{amsmath}
\usepackage{amssymb, amsfonts,amscd,verbatim}

\usepackage[normalem]{ulem}
\usepackage{hyperref}
\usepackage{indentfirst}
\usepackage{latexsym}
\input xy
\xyoption{all}

\usepackage{amsmath}    

\theoremstyle{plain}
\newtheorem{Pocz}{Poczatek}[section]
\newtheorem{Proposition}[Pocz]{Proposition}
\newtheorem{Theorem}[Pocz]{Theorem}
\newtheorem{Corollary}[Pocz]{Corollary}

\newtheorem{Lemma}[Pocz]{Lemma}

\theoremstyle{definition}
\newtheorem{Definition}[Pocz]{Definition}

\theoremstyle{remark}
\newtheorem{Remark}[Pocz]{Remark}

\DeclareMathOperator*{\dist}{dist}
\DeclareMathOperator*{\diam}{diam}

\def\NN{{\mathbb N}}

\def\UU{{\mathcal U}}
\def\VV{{\mathcal V}}
\def\WW{{\mathcal W}}
\def\asdim{\mathrm{asdim}}

\def\dim{\mathrm{dim}}
\def\diam{\mathrm{diam}}
\def\dokaz{{\bf Proof. }}
\def\edokaz{\hfill $\blacksquare$}

\numberwithin{equation}{section}

\title[
Large scale absolute extensors
]%
  {Large scale absolute extensors}

\author[J. Dydak]{Jerzy Dydak}
\address{Department of Mathematics\\227 Ayres Hall\\ University of Tennessee\\ Knoxville, TN 37996,USA}
\email{jdydak@utk.edu}

\author[A.J. Mitra]{Atish J. Mitra}
\address{Department of Mathematical Sciences\\ Montana Tech of The University of Montana\\ 1300 W. Park Street\\
Butte, MT 59701, USA}
\email{atish.mitra@gmail.com}

\date{ \today
}


\begin{document}


\baselineskip=17pt


\maketitle

\begin{abstract}
This paper is devoted to dualization of dimension-theoretical results from the small scale to the large scale.
So far there are two approaches for such dualization: one consisting of creating analogs of small scale concepts and the other amounting to the covering dimension of the Higson corona $\nu(X)$ of $X$. The first approach was used by M.Gromov when defining 
the asymptotic dimension $\asdim(X)$ of metric spaces $X$. The second approach was implicitly contained in the paper 
\cite{Dran AsyTop} by Dranishnikov on asymptotic topology. It is not known if the two approaches yield the same concept. However, Dranishnikov-Keesling-Uspenskiy proved $\dim(\nu(X)\leq \asdim(X)$ and Dranishnikov established that $\dim(\nu(X)= \asdim(X)$ provided $\asdim(X) < \infty$. We characterize asymptotic dimension (for spaces of finite asymptotic dimension) in terms of extensions of slowly oscillating functions to spheres. Our approach is specifically designed to relate asymptotic dimension to the covering dimension of the Higson corona $\nu(X)$ in case of proper metric spaces $X$.  As an application, we recover the results of Dranishnikov-Keesling-Uspenskiy 
and Dranishnikov.
\end{abstract}

\section{Introduction}

\textbf{Asymptotic dimension} of metric spaces was introduced by M. Gromov \cite{Grom} as a means of exploring large scale properties of the space and has been studied extensively during the last two decades. Gromov's definition (see \ref{asdim via UBDcovers with arb high Lebesgue}) dualizes covering dimension as follows: instead of refining open covers by open covers of multiplicity at most $n+1$,
it asks for coarsening of uniformly bounded covers $\mathcal{U}$ by uniformly bounded covers $\mathcal{V}$ with the property that
every element $U$ of $\mathcal{U}$ intersects at most $n+1$ elements of $\mathcal{V}$.

In case of proper metric spaces $X$ (that means bounded subsets of $X$ have compact closure) there is another way to introduce a coarse invariant related to dimension. Namely, it is the covering dimension $\dim(\nu(X))$ of the Higson corona $\nu(X)$ of $X$ (see below). As shown in \cite{Roe lectures}, two proper metric spaces that are coarsely equivalent have homeomorphic Higson coronas, hence their covering dimensions are the same. Thus, indeed, $\dim(\nu(X))$ is a coarse invariant and it is an open question if it is equal to the asymptotic dimension of $X$. This paper is devoted to the internal characterization of $\dim(\nu(X))$, a characterization that makes sense to all metric spaces, not just proper metric spaces.

We know (see \cite{DranKeesUsp}) that  for a proper metric space $X$, the covering dimension of the Higson corona does not exceed the asymptotic dimension of $X$. Also, see \cite{Dran AsyTop}, $\dim(\nu(X)=\asdim(X)$ if $X$ is a proper metric spaces of finite asymptotic dimension $\asdim(X)$. Our approach gives alternative proofs of those results.

Recall that the Higson corona $\nu(X)$ is the complement of $X$ in its Higson compactification $h(X)$. In turn, the Higson compactification $h(X)$ is characterized by the fact that all continuous slowly oscillating functions from $X$ to the unit interval $[0,1]$ extend to continuous functions on $h(X)$. Thus, it is an analog of the \v Cech-Stone compactification where the family of all continuous functions is replaced by all continuous and slowly oscillating functions. The simplest definition of $f:X\to [0,1]$ being \textbf{slowly oscillating} is that $|f(x_n)-f(y_n)|\to 0$ whenever $\sup(d(x_n,y_n)) < \infty$ and $x_n\to \infty$ (that means each bounded subset of $X$ contains only finitely many elements of the sequence $\{x_n\}_{n\ge 1}$). In particular, every function of compact support is slowly oscillating, hence the Higson compactification of $X$ does exist and contains $X$ topologically.

Since covering dimension $\dim(X)$ of a compact space $X$ being at most $n$ can be characterized by saying that the $n$-sphere $S^n$ is an absolute extensor of $X$ (that means any continuous map $f:A\to S^n$, $A$ a closed subset of $X$, can be extended over $X$), one should look for analogous concept involving slowly oscillating functions (but not necessarily continuous). Therefore we introduce the concept of a \textbf{large scale absolute extensor} of a metric space. $K$ is a {\bf large scale absolute extensor of $X$} ($K\in \textrm{ls-AE}(X)$) if for any subset $A$ of $X$ and any slowly oscillating function $f\colon A\to K$ there is a slowly oscillating extension $g\colon X\to K$. We characterize large scale absolute extensors of a space in terms of extensions of $(\epsilon,R)$-continuous functions. It turns out that being large scale absolute extensor of  a space is a coarse invariant of the space. In the later part of the paper we find necessary and sufficient conditions for a sphere $S^m$ to be a large scale extensor of $X$. This is done by comparing existence of Lipschitz extensions in a finite range of Lipschitz constants to existence of Lebesgue refinements in a finite range of Lebesgue constants. We characterize asymptotic dimension of the space in terms of spheres being large scale absolute extensors of the space.

Another natural idea is to study large scale extensors $Y$ of a metric space $X$ defined as follows:

for any $\epsilon > 0$ there is $\delta > 0$ such that any $(\delta,\delta)$-Lipschitz map $f\colon A\subset X\to Y$ extends to a $(\epsilon,\epsilon)$-Lipschitz map $g\colon X\to Y$. It turns out, for bounded metric spaces $Y$, this approach is equivalent to the one involving slowly oscillating functins.

There have been two earlier approaches where connections of asymptotic dimension with extensions of other categories of maps have been studied. One is A.Dranishnikov's (\cite{Dran AsyTop}, using extensions of proper asymptotically Lipschitz functions to euclidean spaces), and by Repov\v s-Zarichnyi (\cite{RepovsZarichnyi_OnAsymptoticExtensionDimension} using maps to open cones). In a separate paper \cite{DydakMitra lsAE-2} we examine how those existing concepts relate to our concept of large scale absolute extensors. 

\subsection{C*-algebra approach to asymptotic dimension}

There is another way to create large scale analogs of covering dimension by using the concept of the \textbf{nuclear dimension} of C*-algebras (see \cite{WinZach}). Given any functor $F$ from the large scale category to the category of C*-algebras, one can consider the nuclear dimension of $F(X)$ for any coarse space $X$ and it gives rise to a large scale invariant. So far we know of two useful functors $F$: one is $\frac{B_h(X)}{B_0(X)}$ (see p.31 of \cite{Roe lectures}), where $B_h(X)$ is the algebra of all bounded functions $X\to \mathbb{C}$ that are slowly oscillating and $B_0(X)$ is the subalgebra of $B_h(X)$ consisting of functions that tend to $0$ at infinity. The second useful C*-algebra of a coarse space $X$ is its uniform Roe algebra (see Chapter 4 of \cite{Roe lectures}).

The algebra $\frac{B_h(X)}{B_0(X)}$ being unital and commutative, has a compact spectrum whose covering dimension coincides with the nuclear dimension of $\frac{B_h(X)}{B_0(X)}$ (see \cite{WinZach}). In case of proper metric spaces $X$, the spectrum of $\frac{B_h(X)}{B_0(X)}$
is exactly the Higson corona $\nu(X)$ of $X$. Thus, one can view the spectrum of $\frac{B_h(X)}{B_0(X)}$ as an abstract Higson corona of any coarse space $X$. A way to interpret our  paper from the point of view of C*-algebras is as computing the nuclear dimension of $\frac{B_h(X)}{B_0(X)}$ via extension properties of slowly oscillating functions.
\section{Basic concepts}

\subsection{The coarse category}
In this section we review basic concepts of the large scale category that are needed for the remaining part of the paper.

Let us modify slightly the concept of a function being bornologous from \cite{Roe lectures}:

\begin{Definition}
Given $\alpha\colon [0,\infty)\to [0,\infty)$ and a function $f\colon (X,d_{X})\to (Y,d_{Y})$ of metric spaces we say $f$ is {\bf $\alpha$-Lipschitz} if $d_{Y}(f(x),f(x'))\leq \alpha(d_{X}(x,x'))$ for any $x,x'\in X$.  

$f\colon (X,d_{X})\to (Y,d_{Y})$ is called {\bf bornologous} if it is $\alpha$-Lipschitz for some non-decreasing $\alpha\colon [0,\infty)\to [0,\infty)$.
\end{Definition}

Notice $\lambda$-Lipschitz functions correspond to $\alpha$ being the dilation $\alpha(x)=\lambda\cdot x$
 and $(\lambda,C)$-Lipschitz functions (or \textbf{asymptotically Lipschitz} functions) correspond to $\alpha(x)=\lambda\cdot x+C$.

We identify two bornologous functions
$f\colon X\to Y$ that are within finite distance from each other and that leads to the
\textbf{coarse category} (or the \textbf{large scale category}).

A function $f\colon X\to Y$ is a \textbf{coarse embedding} (or a \textbf{large-scale embedding}) if there are non-decreasing functions $\alpha\colon [0,\infty)\to [0,\infty)$ and $\beta\colon [0,\infty)\to [0,\infty)$  with $\displaystyle{\lim_{t \to \infty} \alpha(t) = \infty}$ such that
$$\alpha(d_X(x,x'))\leq d_Y(f(x),f(x'))\leq \beta(d_X(x,x'))$$
for all $x,x'\in X$.

A function $f\colon X\to Y$ is a \textbf{ coarse equivalence} (or a \textbf{coarse isomorphism}  if it is a coarse embedding and there is a constant $D >0$ such that any point of $f(X)$ is $D$-close to a point of $Y$.

\subsection{Large scale Continuity}
The aim of this section is to dualize the concept of continuity to the large scale.

\begin{Definition} \cite{CDV3}
Let $\epsilon, \delta >0$. A function $f:(X,d_X)\to (Y,d_Y)$ of metric spaces is $(\epsilon,\delta)$-\textbf{continuous}
if $d_X(x,y) < \delta$ implies $d_Y(f(x),f(y)) < \epsilon$ for all $x,y\in X$.
\end{Definition}

\begin{Remark}
The concept of $(\epsilon,\delta)$-continuity coincides with the concept
of $f$ having $(\delta,\epsilon)$-variation (see \cite{Willett}).
\end{Remark}

The following is a dualization of the standard definition of uniformly continuous functions:

\begin{Definition}\label{lscontinuous}
A function $f:(X,d_X)\to (Y,d_Y)$ of metric spaces is \textbf{large scale continuous} if and only if for every $\delta>0$ there is $\epsilon>0$ such that $f$ is $(\epsilon,\delta)$-continuous.
\end{Definition}

It turns out large scale continuity of $f$ is equivalent to $f$ being bornologous. However, large scale continuity, in addition to being dual to uniform continuity, is easier to apply in some cases.

\begin{Proposition}
 $f:(X,d_X)\to (Y,d_Y)$ of metric spaces is large scale continuous if and only if 
 it is bornologous.
\end{Proposition}
\begin{proof}
 In one direction the proof is immediate: if $f$ is $\alpha$-Lipschitz,
 then $\epsilon=\alpha(\delta)$ works. 
 
 Assume $f$ is large scale continuous and define $\alpha(\delta)$ as follows:\\
 1. $\alpha(0)=0$.\\
 2. if $\delta > 0$, then $\alpha(\delta)$ is the infimum of all $\epsilon > 0$
 such that $f$ is $(\epsilon,\delta)$-continuous.
 
 Notice $\alpha$ is non-decreasing and $f$ is $\alpha$-Lipschitz.
\end{proof}

\subsection{Slowly oscillating functions}

Let us generalize the usual concept of slowly oscillating functions to functions between any metric spaces. 

\begin{Definition}
Let $X,Y$ be metric spaces and $x_0 \in X$.  We say that a function $f: X \to Y$ is \textbf{slowly oscillating} if for any $R,\epsilon > 0$ there is a $N>0$ such that for any $x \in X$ with $d(x_0,x)>N$  the diameter of the set $f(B(x_0,R))$ is less than $\epsilon$. 
\end{Definition}

We note that in the above definition of slowly oscillating functions we do not require the  functions to be continuous.

\begin{Proposition}
Two metrics $d_X$ and $\rho_X$ on $X$ are large scale equivalent if and only if they have the same bounded sets and
any function $f\colon X\to K$ that is slowly oscillating with respect to one metric is also slowly oscillating with respect
to the other metric.
\end{Proposition}
\dokaz
Assume $d_X$ and $\rho_X$ on $X$ have the same bounded sets.
Notice $d_X$ and $\rho_X$ on $X$ are large scale equivalent if and only if
any family of sets that is uniformly bounded with respect to one metric is also uniformly bounded with respect to the other metric. It follows that any function $f\colon X\to K$ that is slowly oscillating with respect to one metric is also slowly oscillating with respect
to the other metric.

Assume any function $f\colon X\to K$ that is slowly oscillating with respect to one metric is also slowly oscillating with respect
to the other metric. Suppose $d_X$ and $\rho_X$ are not large scale equivalent. Without loss of generality
assume there is a sequence $\{(x_n,y_n)\}$ in $X\times X$ such that $\{d_X(x_n,y_n)\}$ is bounded,
$\rho_X(x_n,y_n)\to\infty$ and $x_n\to\infty$ with respect to metric $d_X$.
The function $f\colon X\to [0,1]$ sending all $x_n$'s to $0$ and sending all $y_n$'s to $1$ is slowly oscillating with respect to
$\rho_X$, so it can be extended to a slowly oscillating function $F\colon (X,\rho_X)\to [0,1]$.
Notice $F$ is not slowly oscillating with respect to $d_X$, a contradiction.
\edokaz

\begin{Proposition}
Two metrics $d_X$ and $\rho_X$ on $X$ are uniformly equivalent if and only if
any function $f\colon K\to X$ that is slowly oscillating with respect to one metric is also slowly oscillating with respect
to the other metric.
\end{Proposition}
\dokaz
If metrics $d_X$ and $\rho_X$ on $X$ are not uniformly equivalent, then there is $\epsilon > 0$
such that distances from $x_n$ to $y_n$ with respect to one metric are converging to $0$ yet
distances from $x_n$ to $y_n$ with respect to the other metric are all greater than $\epsilon$.
Put $K=\{0,1\}\times \{n^2\}_{n=1}^\infty$ and define $f\colon K\to X$
by $f(0,n^2)=x_n$, $f(1,n^2)=y_n$ for all $n\ge 1$. Notice $f$ is slowly oscillating with respect to one metric
only.
\edokaz

\begin{Definition}
Given a metric space $(K,d_K)$ its $M$-\textbf{micro-version} is $(K,d_K^M)$,
where $d_K^M(x,y)=d_K(x,y)$ if $d_K(x,y)\leq M$ and $d_K^M(x,y)=M$ if $d_K(x,y)\ge M$.
\end{Definition}

\begin{Definition}
Given a metric space $(X,d_X)$ its $M$-\textbf{macro-version} is $(X,d_X^M)$,
where $d_X^M(x,y)=d_X(x,y)$ if $d_K(x,y)\ge M$ and $d_X^M(x,y)=0$ if $0 < d_K(x,y)< M$.
\end{Definition}

\subsection{Higson Compactification}

We will use the following characterization of the Higson compactification from \cite{Kees}. For definition and details about the Higson compactification see \cite{Kees}.

\begin{Proposition}\label{Keesling Characof Higson Compactification}
Suppose that X is a noncompact proper metric space. The Higson compactification  $h(X)$
is the unique compactification of $X$ such that if $Y$ is any compact
metric space and $f: X \to Y$ is continuous, then $f$ has a
continuous extension to $h(X)=X\cup \nu(X)$ if and only if $f$ is slowly oscillating.
\end{Proposition}

\subsection{Concepts related to covers}

We introduce the concept of \emph{dimension of a cover} of a set.
Traditionally, one mingles multiplicity and dimension. One might as well extend the concept of dimension from spaces to covers to simplify exposition of results and proofs.
\begin{Definition}\label{DimCoverDef}
If $\mathcal{U}$ is a family of subsets of a set $X$, then $\dim(\mathcal{U})\leq n$ means that each $x\in X$ is contained in at most $(n+1)$ elements of $\mathcal{U} $. Equivalently, the \textbf{multiplicity} $m(\mathcal{U} )$ of $\mathcal{U} $ is at most $n+1$.
\end{Definition}

\begin{Definition}\label{LebesgueDef}
The \textbf{Lebesgue number} $Leb(\mathcal{U}) $ of a cover $\mathcal{U} $  of
$X$ is the supremum of all $r \ge 0$ such that every $r$-ball $B(x,r)$ is contained in some element of $\mathcal{U} $.
\end{Definition}

If the Lebesgue number of $\mathcal{U} $ is at least $R$, we express it by saying $\mathcal{U} $ is $R$-\textbf{Lebesgue}.

\begin{Definition}
The \textbf{diameter} $\diam(\mathcal{U})$ of a family of sets in a metric space $X$ is the supremum of distances $d(x,y)$, where $x$ and $y$ belong to the same element of $\mathcal{U}$.\\
$\mathcal{U} $ is \textbf{uniformly bounded} (or $M$-\textbf{bounded}) if $\diam(\mathcal{U} ) < \infty$ (if $\diam(\mathcal{U} ) < M$).
\end{Definition}

\section{Dualizing small scale dimensions}

There are three major ways to define covering dimension for normal topological spaces $X$:\\

1. In terms of open covers of $X$;

2. In terms of pushing maps $f:X\to K$ (from $X$ to a CW complex $K$) into the $n$-skeleton $K^{(n)}$ of $K$;

3. In terms of extending maps $f:A\to S^n$ (from a closed subset $A$ of $X$ to the $n$-sphere) over the whole $X$.\\

One way of dualizing 1) to the coarse category was given by \cite{Grom} (see \ref{asdim via UBDcovers with arb high Lebesgue} below). 2) was dualized in  \cite{CDV2} and \cite{Dran AsyTop} (see also \cite{RepovsZarichnyi_OnAsymptoticExtensionDimension}) has a dualization of 3) involving $R^{n+1}$ instead of $S^n$.

This paper is devoted to a different way of dualizing 3), one that uses extensions of slowly oscillating functions instead of continuous maps. It turns out that generalization is related to an alternative way of dualizing 1).

There are two ways of defining covering dimension of a space $X$ using covers:

a. any open cover of $X$ admits an open refinement of dimension at most $n$,

b. any finite open cover of $X$ admits an open refinement of dimension at most $n$.

Definition b) works well in case of normal spaces $X$ and is equivalent to $S^{n}$ being an absolute extensor of $X$. Definition a) works well in case of paracompact spaces $X$ and is equivalent to b) in that case.

We consider the following possible dualizations of the above two cases to the large scale category.

A.  For every $r > 0$ there is $s > 0$ such that any cover $\UU$ of $X$ of Lebesgue number at least $s$ admits a refinement $\VV$ of Lebesgue number at least $r$ and dimension at most $n$.

B.  For every $r > 0$ there is $s > 0$ such that any finite cover $\UU$ of $X$ of Lebesgue number at least $s$ admits a refinement $\VV$ of Lebesgue number at least $r$ and dimension at most $n$.

To compare A) to Gromov's original definition of asymptotic dimension, let us recall one of the many equivalent definitions of asymptotic dimension of a metric space:

\begin{Definition}\cite{BD1}\label{asdim via UBDcovers with arb high Lebesgue}

A metric space $X$  is of \textbf{asymptotic dimension} at most $n$ if and only if for every $\lambda >0$ there exists a uniformly bounded cover of $X$ with Lebesgue number  at least $\lambda$ and dimension at most $n$.

\end{Definition}

The following proposition shows that A) makes a good definition for asymptotic dimension for any metric space.

\begin{Proposition}\label{Dualization1 for asdim}

For a metric space $X$,   $\asdim X \le n$

if and only for every $R > 0$ there is $S > 0$ such that  any $S$-Lebesgue cover of $X$ admits a $R$-Lebesgue refinement of dimension at most $n$.

\end{Proposition}

\dokaz

If $\asdim X \le n$ then for any $R>0$ using the definition \ref{asdim via UBDcovers with arb high Lebesgue} we get an uniformly bounded $R$-Lebesgue cover $\UU$ of dimension at most $n$. Any  $S$-Lebesgue cover $\VV$ with  $S > \text{mesh}(\UU)$ will have $\UU$ as a refinement.

Conversely, for any $R > 0$ get a $S>0$ using hypothesis. Then the cover of $X$ by $S$-balls around each point has an $R$-Lebesgue refinement of dimension at most $n$.

\edokaz

We will see in a later section (see \ref{Dualization2forasdim}) the closest we can get to definition B).

It is worth comparing the above characterization of asymptotic dimension to another dimension from the small scale world. In  \cite{Isbell}, the author defines uniform dimensions of uniform spaces $X$. For metric spaces this takes two forms. Following \cite{Isbell} we define uniform coverings to be covers with positive Lebesgue numbers.

\begin{Definition}\cite{Isbell}\label{uniform dimensions}

\begin{enumerate}

\item A metric space   is of large uniform dimension at most $n$ ($\Delta d (X) \le n$) if and only if for any  uniform covering  there exists a uniform refinement  with dimension at most $n$.

\item A metric space   is of (small) uniform dimension at most $n$ ($\delta d (X) \le n$) if and only if for any  finite uniform covering  there exists a uniform refinement  with dimension at most $n$.

\end{enumerate}

\end{Definition}

In the next proposition we note that the definition of large uniform dimension is a direct (but trivial)  dual of our characterization of asymptotic dimension 3.2.

\begin{Proposition}\label{Dualization for large unif dimension }

For a metric space $X$,   $\Delta d (X) \le n$ iff for every $S > 0$ there is $R > 0$ such that  any $S$-Lebesgue cover of $X$ admits a $R$-Lebesgue refinement of dimension at most $n$.

\end{Proposition}

\section{Large scale absolute extensors}
In this section we create analogs in the large scale category of absolute extensors in the topological category. 

\begin{Definition}
A metric space $K$ is a {\bf large scale absolute extensor} of a metric space $X$ (notation: $K\in \textrm{ls-AE}(X)$) if for any subset $A$ of $X$
and any slowly oscillating function $f\colon A\to K$ there is an extension $g\colon X\to K$
of $f$ that is  slowly oscillating.
\end{Definition}

\textbf{Warning}: Functions in the above definition are not assumed to be continuous.

\begin{Proposition}
The real line $R$ is not a large scale absolute extensor of itself.
\end{Proposition}
\dokaz
Let $A$ be the subset of $R$ consisting of squares of all integers. The inclusion $i\colon A\to R$
is slowly oscillating as any sequence $(x_n,y_n)\in A\times A$ diverging to infinity such that $\{|x_n-y_n|\}$ is bounded
must be on the diagonal of $A\times A$ starting from some $n$.
Suppose $i$ extends to a slowly oscillating function $f\colon R\to R$.
There is $M > 0$ such that $|f(n+1)-f(n)| < \frac{1}{2}$ for $n > M$.
Therefore,
$$(n+1)^2-n^2=|f((n+1)^2)-f(n^2)| \leq \sum\limits_{i=n^2}^{i=(n+1)^2-1} |f(i+1)-f(i)|< \frac{(n+1)^2-n^2}{2}$$ for $n > M$, a contradiction.
\edokaz

The following result shows that one may restrict attention to complete metric spaces $K$
when discussing large scale absolute extensors.

\begin{Proposition}
If $L$ is dense in $K$, then the following conditions are equivalent for any metric space $X$:
\begin{itemize}
\item[a.] $K$ is a large scale extensor of $X$,
\item[b.] $L$ is a large scale extensor of $X$.
\end{itemize}
\end{Proposition}
\dokaz
Suppose $f\colon (X,A)\to (K,L)$.
By an approximation
$g\colon X \to L$ of $f$ at infinity we mean a function $g$ such that $g|A=f|A$ and $d(f(x),g(x))\to 0$ as $x\to \infty$.
To construct $g$
pick $x_0\in X$ and for each $x\in X\setminus A$ pick $g(x)\in L$ such that $d_K(g(x),f(x)) < \frac{1}{1+d_X(x,x_0)}$.

a)$\implies$b). Suppose $f\colon A\subset X\to L$ is slowly oscillating and choose an extension
$g\colon X\to K$ of $f$ that is slowly oscillating. Choose an approximation $h\colon X\to L$ of $g$
at infinity such that $h|A=g|A$. Notice $h$ is a slowly oscillating extension of $f$.
\par b)$\implies$a). Suppose $f\colon A\subset X\to K$ is slowly oscillating and choose
an approximation $g\colon A\to L$ of $f$
at infinity. Notice $g$ is a slowly oscillating, so it has a slowly oscillating extension of $h\colon X\to L$.
Paste $f$ and $h|(X\setminus A)$ to obtain a slowly oscillating extension $F\colon X\to K$ of $f$.
\edokaz

\begin{Corollary}
Being a large scale absolute extensor of a metric space $X$ is an invariant in the uniform
category.  Being a compact large scale absolute extensor of a metric space $X$ is an invariant in the topological
category.
\end{Corollary}

Our next result shows one can reduce investigation of large scale extensors to bounded metric spaces $K$
and discrete metric spaces $X$.

\begin{Corollary}\label{MacroMicroReduction}
Given metric spaces $X$ and $K$ the following conditions are equivalent:
\begin{itemize}
\item[a.] $K$ is a large scale extensor of $X$.
\item[b.] Any micro-version of $K$ is a large scale extensor of $X$.
\item[c.] $K$ is a large scale extensor of any macro-version of $X$.
\end{itemize}
\end{Corollary}

The following result gives a characterization of large scale absolute extensors in terms of the extensions of (-,-)-continuous functions.
Its importance lies in the fact that it relates $K$ being a large scale extensor of $X$ to the behavior of
all bounded subsets of $X$.

\begin{Theorem}\label{BigTheorem}
The following conditions are equivalent:
\begin{itemize}
\item[a.] $K\in \textrm{ls-AE}(X)$,
\item[b.] For all $M,\epsilon > 0$ there is $n, R, \delta > 0$ such that for any bounded subset $B$ of
$X\setminus B(x_{0},n)$ any $(\delta,R)$-continuous function $f\colon A\subset B\to K$
extends to an $(\epsilon,M)$-continuous function $g\colon B\to K$.
\item[c.] For all $M,\epsilon > 0$ there is $R, \delta > 0$ such that for any bounded subset $B$ of
$X$ any $(\delta,R)$-continuous function $f\colon A\subset B\to K$
extends to an $(\epsilon,M)$-continuous function $g\colon B\to K$.
\item[d.] For all $M,\epsilon > 0$ there is $R, \delta > 0$ such that
any $(\delta,R)$-continuous function $f\colon A\subset X\to K$
extends to an $(\epsilon,M)$-continuous function $g\colon X\to K$.
\end{itemize}
\end{Theorem}
\dokaz
a)$\implies$b). Suppose there is $M,\epsilon > 0$ with the property that for any choice of $n, R, \delta > 0$
there is a bounded subset $B$ of
$X\setminus B(x_{0},n)$ and an $(\delta,R)$-continuous function $f\colon A\subset B\to K$
but no extension $g\colon B\to K$ of $f$ is $(\epsilon,M)$-continuous.

By induction, as described below, choose a sequence of functions $f_{n}\colon A_{n}\subset B_{n}\to K$ that are $(\frac{1}{n},n)$-continuous, do not extend over $B_{n}$ to an $(\epsilon,M)$-continuous function,
and the distance between any two points $x\in B_{i}$, $y\in B_{j}$ is at least $i$ if $i < j$.

First choose bounded subsets $A_1 \subset B_1 \subset X\setminus B(x_{0},1)$ and a $(1,1)$-continuous function $f_{1}\colon A_{1}\subset B_{1}\to K$
that  does not extend over $B_{1}$ to an $(\epsilon,M)$-continuous function. Now suppose functions $\{f_i\}_1^n$ have been chosen satisfying the above conditions, and we want to create $f_{n+1}$. Choose $N$ such that $B_n \subset B(x_0, N) $ and choose bounded subsets $A_{n+1} \subset B_{n+1} \subset X\setminus B(x_{0},N+n)$ and a  function $f_{n+1}\colon A_{n+1}\subset B_{n+1}\to K$
that is $(\frac{1}{n+1},n+1)$-continuous and  does not extend over $B_{n+1}$ to an $(\epsilon,M)$-continuous function.

Paste all $f_{n}$ to $f\colon A=\bigcup\limits_{n=1}^{\infty} A_{n}\to K$ and notice $f$ is slowly oscillating.
Therefore it extends to a slowly oscillating $g\colon X\to K$.
Since each $g\vert B_{n}$ is not $(\epsilon,M)$-continuous, there are points $x_{n},y_{n}\in B_{n}$ for each $n\ge 1$
such that $d_{X}(x_{n},y_{n})\leq M$ but $d_{K}(f(x_{n}),f(y_{n})) > \epsilon$. That contradicts
$g$ being slowly oscillating.
\par
b)$\implies$c). Suppose $M,\epsilon > 0$. Choose $n, S, \mu > 0$ such that
for any bounded subset $B$ of
$X\setminus B(x_{0},n)$ any function $f\colon A\subset B\to K$ that is $(\mu,S)$-continuous
extends to an $(\epsilon,M)$-continuous function $g\colon B\to K$.

Choose $m, T, \lambda > 0$ such that
for any bounded subset $B$ of
$X\setminus B(x_{0},m)$ any function $f\colon A\subset B\to K$ that is $(\lambda,T)$-continuous
extends to a function $g\colon B\to K$ that is $(\mu,S)$-continuous. We may increase $T$ and $S$, so assume $T > S > M$.

Put $\delta =\min(\mu/2,\lambda)$ and put $R=m+3T$. Assume $f\colon A\subset B\to K$ is $(\delta,4R)$-continuous
and $B$ is bounded.

{\bf Case 1}: $A\cap B(x_{0},R)=\emptyset$. Extend $f$ over $A\cup B(x_{0},m+2T)\setminus B(x_{0},m)\to K$
by sending $B(x_{0},m+2S)\setminus B(x_{0},m)$ to a set of diameter $0$.
Notice any such extension $f_{1}$ is $(\lambda,T)$-continuous. Extend $f_{1}$ to $h\colon B \cup B(x_{0},m+2T)\setminus B(x_{0},m)\to K$
that is $(\mu,S)$-continuous. Extend $h$ over $B(x_{0},m+2T)$ by requiring that set is sent to a subset of diameter $0$.
Notice any such extension is $(\mu,S)$-continuous.

{\bf Case 2}: There is a point $x_{1}\in A\cap B(x_{0},R)$.
Notice $\diam(f(A\cap B(x_{0},R)) < \delta$ and let $f_{1}\colon A_{1}=A\cup B(x_{0},m+2S)\to K$
be the extension of $f$ such that $f_{1}(B(x_{0},m+2S)\setminus A)\subset \{f(x_{1})\}$.
Observe $f_{1}$ is $(\mu,S)$-continuous. Indeed, the most relevant case is that of
points $y\in A\setminus B(x_{0},m+2S)$ and $x\in B(x_{0},m+2S)$ such that
$d_{X}(x,y)\leq S$. In that case $d_{X}(y,x_{1}) < R+m+2S+S < 4R$,
so $d_{K}(f_{1}(y),f_{1}(x)) < 2\delta \leq \mu$.
Extend $f_{1}\vert A_{1}\setminus B(x_{0},m)$ to $g_{1}\colon B\cup B(x_{0},m+2S)\setminus B(x_{0},m)\to K$
so that $g_{1}$ has $(M,\epsilon)$-continuous. Pasting $g_{1}$ with $f_{1}$ gives an extension
of $f$ over $B$ that is $(\epsilon,M)$-continuous.
\par
c)$\implies$d). Suppose $M,\epsilon > 0$. Choose $S, \mu > 0$ such that
for any bounded subset $B$ of
$X$ any function $f\colon A\subset B\to K$ that is $(\mu,S)$-continuous
extends to a function $g\colon B\to K$ that is $(\epsilon,M)$-continuous.

Choose $T, \lambda > 0$ such that
for any bounded subset $B$ of
$X$ any function $f\colon A\subset B\to K$ that is $(\lambda,T)$-continuous
extends to a function $g\colon B\to K$ that is $(\mu,S)$-continuous. We may increase $T$, so assume $T > S$.

Put $\delta =\min(\mu,\lambda)$ and put $R=3T$. Assume $f\colon A\to K$ is $(\delta,4R)$-continuous.
Put $C_{k}=B(x_{0},(2k+2)R)\setminus B(x_{0},(2k-1)R)$ for $k\ge 0$.
There is an extension $g_{k}\colon C_{k}\to K$ of $f \vert C_{k}\cap A$ that is $(\mu,S)$-continuous.

Paste $g_{k}\vert (B(x_{0},(2k+1)R)\setminus B(x_{0},2kR))$ with $g_{k+1}\vert (B(x_{0},(2k+3)R)\setminus B(x_{0},(2k+2)R))$
and with $f\vert (B(x_{0},(2k+3)R)\setminus B(x_{0},2kR))\cap A$ to obtain a function that is
$(\mu,S)$-continuous, so it extends over $B(x_{0},(2k+3)R)\setminus B(x_{0},2kR)$ to a function
$h_{k}$ that is $(\epsilon,M)$-continuous. Pasting all $h_{k}$ together produces an extension of $f$
that is $(\epsilon,M)$-continuous.

\par
d)$\implies$a). Suppose there is a slowly oscillating function $f:A \subset X \to K$. For every $n $, there is $S_n,\mu_n >0$ such that any function $g:A \subset X \to K $ that is $(\mu_n,S_n)$-continuous extends to an $(\frac{1}{n},n)$-continuous function  $\tilde{g}:X  \to K $, and there is $T_n,\lambda_n >0$ such that  any function $g:A \subset X \to K $ that is $(\lambda_n,T_n)$-continuous extends to a $(\mu_n,S_n)$-continuous function  $\tilde{g}:X  \to K $. We can take $T_n > S_n > n$, and also $\{T_n\}$,$\{S_n\}$ to be increasing sequences. We create the extension of function $f$ in two steps. First, find $R_n > 0$ such that for $x,y \in A$ with $d(x_0,x)>R_n$ and $d(x,y)<T_n$ we have $d(f(x),f(y))< \lambda_n$. There is an extension of  $f_n:A \cap (B(x_0,R_{n+1})\setminus B(x_0,R_{n}))  \to K$ to a $(\mu_n,S_n)$-continuous function $g_n:A \cap (B(x_0,R_{n+1})\setminus B(x_0,R_{n})) \cup (B(x_0,R_{n+1}+n)\setminus B(x_0,R_{n+1}-(n+1)))  \to K$. This defines the first stage function $g$.In the second stage there is an extension of  $g_n:A \cup (B(x_0,R_{n+1}+n)\setminus B(x_0,R_{n+1}-(n+1))) \cup (B(x_0,R_{n-1}+n-1)\setminus B(x_0,R_{n-1}-(n))) \to K$ to a function $h_n:A \cap B(x_0,R_{n+1}+n)\setminus B(x_0,R_{n-1}-n)   \to K$ that is $(\frac{1}{n},n)$-continuous. Finally pasting all the $h_n$ we get the desired slowly oscillating extension $h:X \to K$.
\edokaz

\begin{Lemma}\label{DeltaDeltaMEpsilonVarLemma}
Suppose $\diam(K)\leq M$ and $g\colon Y\to X$ is $\alpha$-Lipschitz, $\alpha:(0,\infty)\to (0,\infty)$.
If $f\colon X\to K$ is $(\delta,\delta)$-Lipschitz, then $f\circ g$ is $(\epsilon,\epsilon)$-Lipschitz
provided $\epsilon < M$ and $\delta < \frac{\epsilon}{\alpha(\frac{M-\epsilon}{\epsilon})+1}$.
\end{Lemma}
\dokaz
We need to show $d_K(f(g(x)),f(g(y))) \leq \epsilon \cdot d_Y(x,y)+\epsilon$ for all $x,y\in Y$.
It is so if $\epsilon \cdot d_Y(x,y)+\epsilon \ge M$, so assume
$\epsilon \cdot d_Y(x,y)+\epsilon < M$ or $d_Y(x,y) <\frac{M-\epsilon}{\epsilon}$.
In this case $d_X(g(x),g(y)) \leq \alpha(\frac{M-\epsilon}{\epsilon})$
and $d_K(f(g(x)),f(g(y))) \leq \delta\cdot \alpha(\frac{M-\epsilon}{\epsilon})+\delta \leq
\delta \cdot (\alpha(\frac{M-\epsilon}{\epsilon})+1)\leq \epsilon\leq \epsilon \cdot d_Y(x,y)+\epsilon$.
\edokaz

\begin{Corollary}
The following conditions are equivalent for a bounded metric space $K$:
\begin{itemize}
\item[a.] $K\in \textrm{ls-AE}(X)$,
\item[b.] For all $\epsilon > 0$ there is $\delta > 0$ such that for any subset $A$ of
$X$ any $(\delta,\delta)$-Lipschitz function $f\colon A\to K$
extends to an $(\epsilon,\epsilon)$-Lipschitz function $g\colon X\to K$.
\end{itemize}

\end{Corollary}
\dokaz Assume $\diam(K) < M$ and $M > 1$.
\par a)$\implies$b). Given $1 > \epsilon > 0$ find $S, \mu > 0$ such that
any $(\mu,S)$-continuous function $f\colon A\to K$  extends to $F\colon X\to K$
that is $(\epsilon,\frac{M-\epsilon}{\epsilon})$-continuous. Put $\delta=\frac{\mu}{S+1}$.
If $f\colon A\to K$ is  $(\delta,\delta)$-Lipschitz, then it is $(\mu,S)$-continuous
as $d_X(x,y)\leq S$ implies $d_K(f(x),f(y))\leq S\cdot \delta + \delta =\delta\cdot(S+1)=\mu$. Pick
an extension $g\colon X \to K$ of $f$ that is $(\epsilon,\frac{M-\epsilon}{\epsilon})$-continuous.
If $d_X(x,y) > \frac{M-\epsilon}{\epsilon}$, then $d_K(g(x),g(y))\leq M\leq \epsilon \cdot d_X(x,y)+\epsilon$.
If $d_X(x,y) \leq \frac{M-\epsilon}{\epsilon}$, then $d_K(g(x),g(y))\leq \epsilon\leq \epsilon \cdot d_X(x,y)+\epsilon$.
\par b)$\implies$a).
Suppose $S, \epsilon > 0$ and put $\mu=\frac{\epsilon}{S+1}$. Pick $1 > \delta > 0$ such that for any subset $A$ of
$X$ any $(\delta,\delta)$-Lipschitz function $f\colon A\to K$
extends to an $(\mu,\mu)$-Lipschitz function $g\colon X\to K$.
Every $(\delta,\frac{M-\delta}{\delta})$-continuous function $f\colon A\subset X\to K$ is $(\delta,\delta)$-Lipschitz,
so it extends to an $(\mu,\mu)$-Lipschitz function $g\colon X\to K$.
Notice $g$ is $(\epsilon,S)$-continuous. By \ref{BigTheorem}, $K\in \textrm{ls-AE}(X)$.
\edokaz

\begin{Corollary}
The unit interval $I=[0,1]$ is a large scale absolute extensor of any metric space $X$.
\end{Corollary}
\dokaz Assume $X$ is $M$-discrete for some $M > 0$ and $f\colon A\subset X\to I$ is $(\delta,\delta)$-Lipschitz.
Notice $f$ is $(\delta+\frac{\delta}{M})$-Lipschitz. By McShane Theorem
(see \cite{McShane} or Theorem 6.2 on p.43 in \cite{Hei})  $f$ extends
to $g\colon X\to I$ that is $(\delta+\frac{\delta}{M})$-Lipschitz.
By choosing $\delta$ sufficiently small, we can accomplish $g$ to be
$(\epsilon,\epsilon)$-Lipschitz.
\edokaz

\begin{Corollary}
The half-open interval $[0,1)$ and the open interval $(0,1)$ are large scale absolute extensors of any metric space $X$.
\end{Corollary}

\section{Spheres as large scale extensors}

The purpose of this section is to find necessary and sufficient conditions
for a sphere $S^m$ to be a large scale extensor of $X$. This is done by comparing existence of Lipschitz extensions  to existence of Lebesgue refinements  (see \ref{ExtnOfLipschitzMaps VS RefinementOfCovers}).

Given a cover $\UU=\{U_s\}_{s\in S}$ of a metric space $(X,d)$
there is a natural family of functions $\{f_s\}_{s\in S}$ associated to $\UU$:
$f_s(x):=\dist(x,X\setminus U_s)$.
If the multiplicity $m(\UU)$ is finite, then $\UU$ has a natural partition of unity
$\{\phi_s\}_{s\in S}$ associated to it:
$$\phi_s(x)=\frac{f_s(x)}{\sum\limits_{t\in S} f_t(x)}.$$
That partition can be considered as a {\bf barycentric map}
$\phi:X\to\NN(\UU)$ from $X$ to the {\bf nerve} of $\UU$. We
consider that nerve with $l_1$-metric. Recall $\NN(\UU)$
is a simplicial complex with vertices belonging to $\UU$
and $\{U_1,\ldots,U_k\}$ is a simplex in $\NN(\UU)$
if and only if $\bigcap\limits_{i=1}^kU_i\ne\emptyset$.

Since each $f_s$ is
$1$-Lipschitz, $\sum\limits_{t\in S} f_t(x)$ is
$2m(\UU)$-Lipschitz and each $\phi_s$ is
$\frac{2m(\UU)}{Leb(\UU)}$-Lipschitz (use the fact that
$\frac{u}{u+v}$ is
$\frac{\max(Lip(u),Lip(v))}{\inf(u+v)}$-Lipschitz). Therefore
$\phi:X\to\NN(\UU)$ is $\frac{4m(\UU)^2}{Leb(\UU)}$-Lipschitz. See
\cite{BD} for more details and better estimates of
Lipschitz constants.

We denote by $\Delta^n$ the standard unit $n$-dimensional simplex
with $l_1$-metric. Its boundary will be
denoted by $S^{n-1}$ or $\partial\Delta^{n}$ (notice that the unit sphere in $R^n$ is
bi-Lipschitz equivalent to the boundary $\partial\Delta^{n}$).

The following will be important in relating large scale extensors to asymptotic dimension.

\begin{Proposition} \label{ExtnOfLipschitzMaps VS RefinementOfCovers}
Suppose $X$ is a metric space,
$m\ge 0$. Then the following are equivalent.

\begin{enumerate}

\item [a.] For any $\epsilon>0$ there is   $\epsilon > \delta > 0$ such that
any $(\delta,\delta)$-Lipschitz function
$f:A\to S^m$, $A$ a subset of $X$, extends to an $(\epsilon,\epsilon)$-Lipschitz function $\tilde f:X\to S^m$.
\item [b.] For any $s>0$ there is $t>s>0$ such that for any finite
$m+2$-element cover $\UU=\{U_0,\ldots,U_{m+1}\}$ of $X$ with Lebesgue number greater than $t$, there is
a refinement $\VV$ so that $\VV$ has Lebesgue number greater than $s$ and the dimension
of $\VV$ is at most $m$.

\end{enumerate}
\end{Proposition}

\dokaz
By switching to a macro-version of $X$ we may assume $X$ is 1-discrete.\\
a)$\implies$b). Let $s>0$, define $\epsilon=\frac{1}{2s(m+1)}$ and get corresponding $\delta$ from hypothesis.  Define $t=\frac{4(m+2)^2}{\delta}$. Considering any $m+2$-element $t$-Lebesgue  cover $\UU$,
get a barycentric map $\phi:X\to \NN(\UU)=\Delta^{m+1}$ with $Lip(\phi)\leq \delta$.
There is $g:X\to \partial\Delta^{m+1}$ such that $Lip(g)\leq 2\epsilon$ and
$g(x)=\phi(x)$ for all $x\in X$ so that $\phi(x)\in \partial\Delta^{m+1}$.
Consider $V_i=\{x\in X\mid g_i(x) > 0\}$.
Notice $\VV=\{V_i\}_{i=0}^{i=m+1}$ is of dimension at most $m$.
Also $x\in V_i$ implies $x\in U_i$, so $\VV$ refines $\UU$.
Given $x\in X$ there is $i$ such that $g_i(x)\ge \frac{1}{m+1}$.
If $d(x,y) < s$, then
$|g_i(x)-g_i(y)| < \frac{1}{m+1}$ and $g_i(y) > 0$. Thus, the ball at $x$ of radius
$s$ is contained in one element of $\VV$.


b)$\implies$a).  For the proof of this direction we think of maps from $X$ to an $(m+1)$-simplex $\Delta^{m+1}$
as a partition of unity. Since we want to create a map to its boundary
$S^m=\partial \Delta^{m+1}$, a geometrical tool is the radial projection $r$
which we splice in the form of $(1-\beta)\cdot r+ \beta\cdot \phi$
with a partition of unity $\phi$ coming from a covering of $X$
of dimension at most $m$.

It follows from \cite{McShane} that there exists $C > 0$
such that given a $\lambda$-Lipschitz $f:A\to \Delta^{m+1}$ one
can extend it to a $C\cdot\lambda$-Lipschitz $g:X\to
\Delta^{m+1}$.

\par Let $\epsilon >0$. Define  $\delta_{1}=\frac{\epsilon}{(m+2)^3(82C+4)}$ and choose $\delta_{2}<\delta_{1}$ such that  for any finite $m+2$-element cover with Lebesgue number greater than $t= \frac{1}{24\delta_{2} C(m+2)}$ there is
a refinement with Lebesgue number greater than $s=\frac{1}{\delta_1}$ and  dimension at most $m$.

We define  $\delta =\min\{\delta_1,\delta_2\}$ and show below that any $(\delta,\delta)$ Lipschitz map $f:A \to S^m$ extends to $(\epsilon, \epsilon)$ $\tilde f:X\to S^m$.

We first extend $f$ to a $2\delta C$ Lipschitz $g:X\to \Delta^{m+1}$.

Let $\alpha:X\to [0,1]$ be defined as $\alpha(x)=(m+2)\cdot \min\{g_i(x)\mid 0\leq i\leq m+1\}$.
Notice $Lip(\alpha)\leq (m+2) 2C\cdot \delta$.
Let $\beta:[0,1]\to [0,1]$ be defined by $\beta(z)=3z-1$ on $[1/3,2/3]$,
$\beta(z)=0$ for $z\leq 1/3$ and $\beta(z)=1$ for $z\ge 2/3$ and note that $Lip(\beta)\leq 3$.

Put $U_i=\{x\in X\mid g_i(x) > \frac{\alpha(x)}{m+2} \text{ or } \alpha(x)> 2/3\}$
and notice $Leb(\UU)\ge r= \frac{1}{24\delta C(m+2)}$ as follows:
\par Case 1: $x\in X$ and $\alpha(x) > 3/4$. Now, for any $y\in X$
with $d(x,y) < \frac{1}{24\delta C (m+2)}$ one has
$\alpha(x)-\alpha(y) \leq 1/12$, which implies that  $\alpha(y) > 2/3$.
Thus, in that case the ball $B(x,\frac{1}{24 \delta C (m+2)})$ is contained in all $U_i$.
\par
Case 2: $\alpha(x)\leq 3/4$. There is $i$ so that $g_i(x)\ge \frac{1}{m+2}$.
Since $\psi_i=g_i-\frac{\alpha}{m+2}$ is $4\delta C$-Lipschitz,
for any $y\in X$ satisfying $d(x,y) < \frac{1}{16\delta C (m+2)}$
one has $\psi_i(x)-\psi_i(y)< \frac{1}{4(m+2)}$
and $\psi_i(y) > 0$ as $\psi_i(x) \ge \frac{1}{4(m+2)}$.

\par

Shrink each $U_i$ to $V_i$ so that $m(\VV)\leq m+1$
and $Leb(\VV)\ge s=\frac{1}{\delta}$.
The barycentric map $\phi:X\to \partial\Delta^{m+1}$
corresponding to $\VV$ has $Lip(\phi)\leq 4(m+2)^2\delta$.

Define $h(x)=\sum\limits_{i=0}^{m+1} (g_i(x)-\frac{\alpha(x)}{m+2})\cdot \frac{1-\beta(\alpha(x))}{1-\alpha(x)}\cdot e_i+\sum\limits_{i=0}^{m+1} \beta(\alpha(x))\cdot \phi_i(x)\cdot e_i$.
To show $Lip(h)\leq \epsilon$
we will use the following observations.
\begin{enumerate}
\item If $u,v:X\to [0,M]$, then $Lip(u\cdot v)\leq M\cdot(Lip(u)+Lip(v))$.
\item In addition, if $v:X\to [k,M]$ and $k > 0$, then $$Lip(\frac{u}{v})\leq M\cdot \frac{Lip(u)+Lip(v)}{k^2}.$$
\item $v(x)=1-\alpha(x)\ge 1/3$ if $\frac{1-\beta(\alpha(x))}{1-\alpha(x)} > 0$.
\end{enumerate}

Therefore $Lip(\sum\limits_{i=0}^{m+1} \beta(\alpha(x))\cdot
\phi_i(x)\cdot e_i)\leq  (m+2)(3Lip(\alpha)+Lip(\phi)) \leq (m+2)(3(m+2)2\delta C + 4(m+2)^2\delta)  \leq (m+2)^3(6C+4)\delta $.

Also,
$Lip(\frac{1-\beta(\alpha(x))}{1-\alpha(x)})\leq 9\cdot 4 \cdot (m+2)\cdot 2\delta \cdot C$, so $Lip(\sum\limits_{i=0}^{m+1}
(g_i(x)-\frac{\alpha(x)}{m+2})\cdot
\frac{1-\beta(\alpha(x))}{1-\alpha(x)})\leq (m+2)\cdot
(4\delta C+72(m+2)\delta C)\leq 76(m+2)^2 \delta C \leq 76(m+2)^3 \delta C $.  So $Lip(h)\leq (m+2)^3(82C+4)\delta \leq \epsilon$.

\par
It remains to show $h(X)\subset \partial\Delta^{m+1}$ and $h|A=f$.
$h|A=f$ follows from the fact $\alpha(x)=0$ if $x\in A$.
It is clear $h(x)\in \partial\Delta^{m+1}$ if either $\beta(\alpha(x))=0$ or $\beta(\alpha(x))=1$,
so assume $0 < \beta(\alpha(x)) < 1$.
In that case $\phi_i(x) > 0$ implies $g_i(x)-\frac{\alpha(x)}{m+2} > 0$,
so the only possibility for $h(x)$ to miss $\partial\Delta^{m+1}$
is when $g_i(x)-\frac{\alpha(x)}{m+2} > 0$ for all $i$ which is not possible.
\hfill $\blacksquare$

From proposition \ref{ExtnOfLipschitzMaps VS RefinementOfCovers} we get the following:

\begin{Corollary} \label{CharOfSmAsLSExtensor}
If $X$ is a metric space and
$m\ge 0$, then the following conditions are equivalent:
\begin{itemize}
\item[a.] $S^m$ is a large scale absolute  extensor of $X$.
\item[b.] For any $s>0$ there is $t > 0$ such that any finite $m+2$-element  cover $\UU=\{U_0,\ldots,U_{m+1}\}$ of $X$
with $Leb(\UU)> t$ admits a refinement $\VV$ so that $Leb(\VV)>s $ and the dimension
of $\VV$ is at most $m$.
\end{itemize}
\end{Corollary}

The next proposition and its corollary will be useful in the proof of \ref{CharacterizingAsymptDimBySpheres}.
\begin{Proposition} \label{nCoversTonPlus1Covers}
Suppose $X$ is a metric space, $n\ge 0$. If for any $s>0$
there is $t>0$ such that every $t$-Lebesgue cover $\UU=\{U_0,\ldots,U_{n+1}\}$ of $X$
admits a $s$-Lebesgue refinement $\VV$ satisfying $m(\VV)\leq n+1$,
then for any $q>0$ there is $r>0$ such that any $r$-Lebesgue cover $\WW=\{W_0,\ldots, W_{n+2}\}$ of $X$
admits a $q$-Lebesgue refinement $\VV$ of dimension at most $n+1$.
\end{Proposition}

\dokaz  By switching to a macro-version of $X$ we may assume $X$ is 1-discrete.\\
Let $q>0$. By hypothesis, there is $t>0$ such that any $n+2$-element $t$-Lebesgue cover $\UU=\{U_0,\ldots,U_{n+1}\}$ of $X$
admits a $q$-Lebesgue refinement $\VV$ satisfying $m(\VV)\leq n+1$.

First we show that any $2t$-Lebesgue $n+2$-element cover
$\UU=\{U_i\}_{i=0}^{i=n+1}$ of $A\subset X$  has a refinement $\VV$ such that
$Leb(\VV)\ge q$ and $m(\VV)\leq n+1$. Define $U'_i=U_i\cup (X\setminus A)$ for $i\leq n+1$ and notice that $Leb(\UU')\ge t$ as follows.
If $x\in X$, then $B(x,t)\cap A$ is either empty
or is contained in $B(y,2t)$ for some $y\in A$.
Since $B(y,2t)\cap A\subset U_i$ for some $i\leq n+1$,
$B(y,2t)\subset U_i'$ and hence $B(x,t)\subset U_i'$. By hypothesis, there is a cover $\WW$ of $X$ such that $\WW$ refines $\UU'$, $Leb(\WW)\ge q$,
and $m(\WW)\leq n+1$. By putting $\VV=\WW|_A$ we get the required refinement.
\par
Suppose $\WW=\{W_0,\ldots,W_{n+2}\}$ is an $r=4t$-Lebesgue cover of $X$.
Let $A$ be the union of balls $B(x,2t)$ such that $B(x,4t)$ is not contained in $W_{n+2}$.
Define $U_i=W_i\cap A$ for $i\leq n+1$ and observe as follows that $Leb(\UU)\ge 2t$
for $\UU=\{U_i\}_{i=0}^{i=n+1}$ as a cover of $A$. If $x\in A$,
then there is $y\in X$ such that $B(y,4t)$ is not contained in $W_{n+2}$
and $x\in B(y,2t)$. Therefore, $B(y,4t)\subset W_i$ for some $i\leq n+1$
which means $B(x,2t)\cap A\subset B(y,4t)\cap A\subset W_i\cap A=U_i$.
\par
Shrink each $U_i$ to $V_i$ so that the intersection of all $V_i$ is empty and
$Leb(\VV)\ge q$.
Define $W'_i=V_i$ for $i\leq n+1$ and $W'_{n+2}=W_{n+2}$.
The cover $\WW'$ is of dimension at most $n+1$.
We show as follows that $Leb(\WW')\ge q$.
If $B(x,4t)\subset W_{n+2}$, we are done. Otherwise $B(x,2t)\subset A$ and
there is $i\leq n+1$ such that $B(x,q)\subset V_i$ in which case $B(x,q)\subset W'_i$.

\hfill $\blacksquare$

\ref{CharOfSmAsLSExtensor} and \ref{nCoversTonPlus1Covers} give the following corollary.

\begin{Corollary} \label{HigherSpheresAreLipAE}
Suppose $X$ is a metric space and $n\ge 0$. If $S^{n}$ is a large scale absolute extensor of $X$, then so is $S^{n+1}$.
\end{Corollary}

\section{Large scale absolute extensors and asymptotic dimension}

The following lemma is a version of the Ostrand-type definition of asymptotic dimension, but with control on Lebesgue number.

\begin{Lemma} \label{OstrandWithLebesgue} A metric space $(X, d)$  is of asymptotic dimension at most $n$ if and only if for every $r >0$ there exist a uniformly bounded covering $\UU$ where $\UU = \bigcup\limits_{i= 1}^{n+1}\UU^i$ with each $\UU^i$ is an $r$-disjoint family and  the Lebesgue number of $\UU$ is at least $r$.
\end{Lemma}

\begin{Proposition} \label{CharacterizingAsymptDimBySpheres}
Suppose $X$ is a metric space of finite asymptotic dimension. If $n\ge 0$, the following conditions are equivalent:

\begin{enumerate}
\item [a.] $S^n$ is a large scale absolute extensor of $X$
\item [b.] $\asdim X \le n$
\end{enumerate}

\end{Proposition}
\dokaz  b)$\implies$a) By taking a 1-net in $X$  we can assume $X$ is 1-discrete.  Let $s>0$ and by hypothesis get a uniformly bounded cover $\VV$ of dimension $\le
n$ with Lebesgue number larger than $s$. If $\UU$ be any cover of Lebesgue number $>t=\text{mesh}(\VV)$. Then $\VV$ is a refinement of $\UU$.

\par a)$\implies$b) In view of \ref{HigherSpheresAreLipAE}, we can assume that $\asdim X \le n+1$. Let $s>0$. There is $t >s> 0$ such that any finite $n+2$-element  cover  of $X$ with $L(\UU)> t$ admits a refinement $\VV$ so that $L(\VV)>s $ and the dimension
of $\VV$ is at most $n$.

As $\asdim X \le n+1$,  there exists an uniformly bounded covering $\UU$ where $\UU = \bigcup_{i= 1}^{n+2}\UU^i$ where each $\UU^i$ is an $t$-disjoint family and  the Lebesgue number of $\UU$ is at least $t$.

Define $U_i$ to be the union of all elements of $\UU^i$ and note that $\{U_i\}_1^{n+2}$ ia $n+2$-element $t$-Lebesgue cover. By hypothesis we can get a $s$-Lebesgue refinement $\VV$ of dimension at most $n$.

For each $V \in \VV$,  there is an $1\le i\le n+2$ such that $V \subset \UU_i$. Replacing $V$ by the collection $\{V\cap W: W \in \UU^i\}$ we get  a $s$-Lebesgue uniformly bounded cover of dimension at most $n$, which implies  $\asdim X \le n$.

\hfill $\blacksquare$

In section 3 we looked for a characterization of asymptotic dimension similar to \ref{Dualization1 for asdim}  using finite covers.

\ref{CharacterizingAsymptDimBySpheres} and \ref{CharOfSmAsLSExtensor} immediately gives  the following:

\begin{Corollary}\label{Dualization2forasdim}
For a metric space $X$ of finite asymptotic dimension,   $\asdim X \le n$
if and only if for every $R > 0$ there is $S > 0$ such that  any finite $n+2$-element $S$-Lebesgue cover of $X$ admits a $R$-Lebesgue refinement of dimension at most $n$.
\end{Corollary}

In the case of covering dimension,  "any finite $n+2$-element cover" can be replaced by "any finite cover". It is not clear if that is true for asymptotic dimension.

\section{Connections to the Higson Corona}

In this section we relate the condition of dimension of the Higson corona of a proper metric space being at most $n$ to $S^n$ being a large scale absolute extensor of $X$.

\begin{Theorem}
If $X$ is a proper metric space and $M$ is a compact ANR, then the following conditions are equivalent.
\begin{enumerate}

\item [a.]  any continuous function $f:A \to M$, where $A$ is a closed subset of $\nu(X)$ extends to a continuous function $F: \nu(X) \to M$

\item [b.] $M$ is a large scale absolute extensor of $X$.
\end{enumerate}
\end{Theorem}

\dokaz
By considering a $1$-net in $X$, we can reduce the proof to a $1$-discrete $X$.

a$\implies$b. Suppose $f:A\subset X\to M$ is slowly oscillating. We will describe how to extend it to a continuous function $F:h(X) \to M$. By the characterizing property of $h(X)$ \ref{Keesling Characof Higson Compactification},  the restriction  $\tilde{f} = F|_X$  is slowly oscillating and is the desired extension of $f$.

 As $M$ is an ANR, we can extend $f$ to a  $g:N\to M$, where $N$ is a closed neighborhood of $A$ in $h(X)$. By hypothesis, we can extend $g|\nu(X)\cap N$ to $G:\nu(X)\to M$, then over a neighborhood $U$ of $\nu(X)$ in agreement with $g$. Now we have to define the extension on $h(X) \setminus U$, which is a compact subset of $X$, hence finite, so we can put any values for the extension there. The resulting function $F:h(X) \to M$ is continuous.

b$\implies$a.

Suppose $f:A\subset \nu(X)\to M$ is continuous. Since $M$ is an ANR, $f$ can be extended  to $g:N\to M$ over a neighborhood  $N$ of $A$ in $h(X)$ .

$g|N\cap X\to M$ is slowly oscillating, so by hypothesis it extends to a slowly oscillating $\tilde{g}:X \to M$. By \ref{Keesling Characof Higson Compactification} $\tilde{g}$  extends to a continuous $G: h(X) \to M$. As $X$ is dense in $h(X)$, $G$ must agree on $A$ with $f$ , so  $G|\nu(X)$ is the required extension of $f$ over $\nu(X)$.
\edokaz

\begin{Corollary}
 If $X$ is a proper metric space and $n\ge 0$, then the following conditions are equivalent:
\begin{enumerate}
\item [a.] $\dim(\nu X)\leq n$,

\item [b.] $S^n$ is a large scale absolute extensor of $X$.
\end{enumerate}
\end{Corollary}

In view of \ref{CharacterizingAsymptDimBySpheres} we have another proof of $\dim(\nu X)=\asdim(X)$
in case of $\asdim(X)$ being finite.

\begin{Corollary}[\cite{DranKeesUsp}  and \cite{Dran AsyTop}]
 If $X$ is a proper metric space and $n\ge 0$, then the following conditions are equivalent:
\begin{enumerate}
\item [a.] $\dim(\nu X)\leq n$,

\item [b.] $\asdim(X)\leq n$.
\end{enumerate}
\end{Corollary}

\subsection*{Acknowledgements}
The second-named author was partially supported by  Israel Science Foundation grant no 836/08.

\end{document}